\documentclass[a4paper,11pt]{article}

\usepackage{authblk}
\usepackage{parskip}
\usepackage{mathtools}
\usepackage{amssymb}
\usepackage{amsmath}
\usepackage{latexsym}
\usepackage{mathrsfs}  % example: \mathscr{A} 
\usepackage{bbm}       % example: \mathbbm{N}
\usepackage{amsthm}    % Theorems
\usepackage{relsize}   % Change default sizes 
\usepackage{graphicx}
\usepackage{thmtools}  % Modify theorem style
\usepackage{mathrsfs}  % Ralph Smith’s Formal Script Font (rsfs)
\usepackage{paralist}  % list on paragraph
\usepackage{lipsum}    % drawing horizontal lines
\usepackage[all]{xy}   % To draw diagrams 

\numberwithin{equation}{section}

\usepackage{titlesec}   %Control what is at the format of section titles

\titleformat{\section}[block]{\bfseries\filcenter}% centred title 
{{\upshape\thesection\enspace}}{.5em}{}

\titleformat{\subsection}[block]{\filcenter}% upshape means that is not in italics
{{\upshape\thesubsection\enspace}}{.5em}{} %%^: These reduce the amount of space taken up by section and subsection

\usepackage{enumitem}  % Modify list
\setlist{nosep}  % nosep= no vertical separation between items.
\setitemize[0]{leftmargin=*}
\setenumerate[0]{leftmargin=*}
\setenumerate[1]{label={(\arabic*)}} % Global setting of enumeration of the form (1), (2), ..
\newcommand{\N}{\mathbb{N}}     % Natural numbers
\newcommand{\R}{\mathbb{R}}     % Real numbers
\newcommand{\Prob}{\mathbb{P}}  % Probability measure
\newcommand{\st}{\,:\,}         % Symbol for "such that"
\newcommand{\inner}[2]{\left\langle #1 \, , \, #2 \right\rangle} % Inner product
\newcommand{\norm}[1]{\left|\left|#1\right|\right|}              % Vector norm
\newcommand{\triplet}[3]{\left( #1, #2, #3 \right) }             % General triplet e.g. a probability space
\newcommand{\ProbSpace}{\triplet{\Omega}{\mathscr{F}}{\Prob}}    % Triplet of a Probability Space
\newcommand{\abs}[1]{\left| #1 \right|}                          % Absolute value  
\renewcommand{\qedsymbol}{$\square$}                       % For a black square at the end of a proof
 
\newcommand{\cadlag}{c\`{a}dl\`{a}g }

\newcommand{\defeq}{\mathrel{\mathop:}=}                         % Defined as equal to symbol
\newcommand\restr[2]{{% we make the whole thing an ordinary symbol
  \left.\kern-\nulldelimiterspace % automatically resize the bar with \right
  #1 % the function
  \vphantom{\big|} % pretend it's a little taller at normal size
  \right|_{#2} % this is the delimiter
  }}
  
\theoremstyle{plain} 
\newtheorem{theo}{Theorem}[section]    

\newtheorem{coro}[theo]{Corollary}
\newtheorem{lemm}[theo]{Lemma}

\theoremstyle{definition} 
\newtheorem{defi}[theo]{Definition}
\newtheorem{exam}[theo]{Example}

\declaretheoremstyle[%
  spaceabove=-5pt,%
  spacebelow=6pt,%
  headfont=\normalfont\itshape,%
  postheadspace=1em,%
  qed=\qedsymbol%
]{mystyle} 
\declaretheorem[name={Proof},style=mystyle,unnumbered,
]{prf}

 \title{Radonification of a cylindrical L\'evy process}

\author{
A. E. Alvarado-Solano\thanks{aalvarado@itcr.ac.cr }   }

\affil{Escuela de Matem\'{a}tica, Instituto Tecnol\'ogico de Costa Rica, 
\\ Cartago, 159-7050, Costa Rica}
\date{}   

\begin{document}

 \maketitle

\abstract{In this work we present a direct proof about radonification of a cylindrical L\'evy process. The radonification technique has been very useful to define an genuine stochastic process starting from a cylindrical process, this is possible thanks to the Hilbert-Schmidt operators. With this work we want to propose a self-contained simple proofs who those who are not familiar with this method, and also present our result which it is apply the radonification method to the case of a cylindrical L\'evy process.}

\smallskip



\section{Introduction}

The radonification technique is the way how some researches had been working to understand the cylindrical noises in linear spaces. This method has been useful in the develop of partial stochastic equations with cylindrical noise, see \cite{LiuRockner}, \cite{RiedleStable-2018}, also it has been used to develop a general theory of integration respect cylindrical martingale-valued measure, see \cite{Alvarado-Fonseca2020} and this technique has been studied in nuclear spaces, see \cite{FonsecaMora:2018-1}.

Our intention on this work is to present results about this radonification method applied to cylindrical L\'evy processes on a Hilbert space. Also we want with this work to let at knowledge to readers about this technique in a self-contained simple proof in the context of a Hilbert space. We consider these proof are valuable because their presentation of the arguments.

This work has three sections. In the second one we present some notation and essential definitions in the following section. In the third and last section we will prove our main result which is the radonification of a cylindrical L\'evy processes by a single Hilbert-Schmidt operator into a genuine L\'evy process. It is known by experts the truthfulness of this fact, however it is unknown to our knowledge any proof of this, but in any case we present a friendly and  original proof, also some consequences of this fact, also we present some typical examples of cylindrical L\'evy processes and how the main result works on them.

\section{Preliminaries}

Let $H$ a Hilbert space and denote by $\langle\cdot,\cdot\rangle$ its inner product and $\Vert\cdot\Vert$ its norm. We will denote by $B_H$ as the unit ball in $H$ and $B_1$ the unit ball in $\mathbb{R}$. We denote by $\mathcal{B}(H)$ its Borel $\sigma$-algebra. We always identify $H$ with its (strong) dual space.  For a probability space $(\Omega, \mathcal{F}, \mathbb{P})$, denote by $L^{0}( \Omega, \mathcal{F}, \mathbb{P}; H)$ (or $L^{0}( \Omega, \mathcal{F}, \mathbb{P})$ if $H=\R$) the linear space of all the equivalence classes of random variables. We equip $L^0(\Omega,\mathcal{F},\mathbb{P};H)$ with the topology of convergence in $\mathbb{P}$-measure and in this case it is complete, metrizable, topological vector space. Similarly, for $p\geq 1$ we denote by $L^p(\Omega,\mathbb{P};H)$ (or $L^p(\Omega,\mathbb{P})$ if $H=\mathbb{R}$) as the space of $p$-integrable random variables which means that $\Vert X\Vert_{L^p(\Omega,\mathbb{P};H)}^p:=\int_{\Omega}{\Vert X\Vert^p}\mathbb{P}(d\omega)<\infty$. With this norm $L^p(\Omega,\mathbb{P};H)$ is a Banach space and for $p=2$ a Hilbert space.

Let $G$ another Hilbert space. The collection of all continuous linear operators will be denoted by $\mathcal{L}(H,G)$. Recall that $S \in \mathcal{L}(H,G)$ is called Hilbert-Schmidt if for some (equivalently for any) complete orthonormal system $(h_{n}: n \in \N)$ in $H$, we have $\Vert S\Vert_{\mathcal{L}_2(H,G)}^2:=\sum_{n=1}^{\infty} \norm{S(h_{n})}^{2}< \infty$. If $S\in\mathcal{L}_2(H,G)$ then there exist complete orthonormal systems $\{h_{n}\}_{ n \in \N}$ in $H$, $\{g_{n}\}_ {n \in \N}$ in $G$, and a sequence $\{\lambda_{n}\}_{n \in \N} \subseteq \R$ such that $\sum_{n=1}^{\infty} \abs{\lambda_{n}}^{2} < \infty$, and:
\begin{equation}\label{spectra}
S(h)= \sum_{n=1}^{\infty} \lambda_{n} \inner{h_{n}}{h} g_{n}, \quad \forall \,  h \in H.
\end{equation}
We will denote by $S^*$ as the adjoint operator of $S$. If $S\in\mathcal{L}_2(H,G)$ then $S^*\in\mathcal{L}_2(G,H)$.

All our random variables will be considered defined on some given complete probability space $\ProbSpace$  equipped with a filtration $(\mathcal{F}_{t} : t \geq 0)$ that satisfies the \emph{usual conditions}, i.e. it is right continuous and $\mathcal{F}_{0}$ contains all subsets of sets of $\mathcal{F}$ of $\Prob$-measure zero. Let $T>0$, we denote $D_T(H)$ as the space of \emph{c\`adl\`ag paths processes} $X=\{X_t\}_{t\in[0,T]}$. This space is equipped with the topology of uniform convergence in probability, and with this it is a complete space. If $t\in\mathbb{R}_+$ we write $D(H)$ for simplify, understand it as a Frech\'et space with the topology generated by the seminorms of $D_{T_n}(H)$ with $T_n\rightarrow\infty$ when $n\rightarrow\infty$.

Let $L=\{L_t\}_{t\in\mathbb{R}_+}$ a $H$-valued processes and we call it a L\'evy process if it initiates at $0 \mathrm{ \ ; \ } \mathbb{P}$-a.s, it is stocastically continuous and has stationary and independents increments. Clasical examples of these are the well knowm Wiener process (or Brownian motion) and the Compound Poisson process. %These processes are well-known by the L\'evy-It\^o decomposition (see \cite{Applebaum-2007}, p.80) which states that  for any $L$ a $H$-valued L\'evy process there exist $b\in H$, a $H$-valued standard  Wiener process $W$ and a Poisson random measure over $\mathbb{R}_+\times H\setminus \{0\}$ such that for each $t\in\mathbb{R}_+$:
%\begin{equation}\label{LevyItoHil}
 %   L_t=tb+W_t+\int_{\Vert h\Vert <1}{h}\tilde{N}_t(dh)+\int_{\Vert h\Vert \geq 1}{h}N_t(dh).
%\end{equation}

%In \eqref{LevyItoHil} each component is independent process from the other, also $\{\int_{\Vert h\Vert<1}{h}\tilde{N}_t(dh)\}_{t\in\mathbb{R}_+}$ is such that $\int_{\Vert h\Vert<1}{h}\tilde{N}_t(dh)=\lim_{n\rightarrow\infty}\int_{\frac{1}{n}<\Vert h\Vert<1}{h}\tilde{N}_t(dh)$ in $L^2(\Omega,\mathbb{P};H)$.

%Also given $T>0$ we denote by $\mathcal{M}_{T}^{2}(H)$ the space of all the $H$-valued c\`adl\`ag square integrable martingales defined on $[0,T]$, which is a Banach space when equipped with the norm (see \cite{DaPratoZabczyk}, Proposition 3.9, p.79), $\norm{M}_{\mathcal{M}_{T}^{2}(H)}^2= \Exp(\sup_{t\in[0,T]} \norm{M_{t}}^{2}) \ $.

For any $n \in \N$ and any $h_{1}, \dots, h_{n} \in H$, we define a linear map $\pi_{h_{1}, \dots, h_{n}}: H \rightarrow \R^{n}$ by $\pi_{h_{1}, \dots, h_{n}}(y)=(\inner{y}{h_{1}}, \dots, \inner{y}{h_{n}})$ $\forall \, y \in H$.  A cylindrical set in $H$ based on $\Gamma \subseteq H$ is a set of the form
$$ \mathcal{Z}\left(h_{1}, \dots, h_{n}; A \right) = \left\{ y \in H \st (\inner{y}{h_{1}}, \dots, \inner{y}{h_{n}}) \in A \right\}= \pi_{h_{1}, \dots, h_{n}}^{-1}(A), $$
where $n \in \N$, $h_{1}, \dots, h_{n} \in \Gamma$ and $A \in \mathcal{B}\left(\R^{n}\right)$. The collection of all the cylindrical sets based on $\Gamma$ is denoted by $\mathcal{Z}(H,\Gamma)$. It is an algebra but if $\Gamma$ is a finite set then it is a $\sigma$-algebra. The $\sigma$-algebra generated by $\mathcal{Z}(H,\Gamma)$ is denoted by $\mathcal{C}(H,\Gamma)$. If $\Gamma=H$, we write $\mathcal{Z}(H)=\mathcal{Z}(H,\Gamma)$ and $\mathcal{C}(H)=\mathcal{C}(H,\Gamma)$. Also if $H$ is separable space we have that $\mathcal{C}(H)=\mathcal{B}(H)$ (see \cite{EvangelioVol1}, Proposition 1.1.1, p.3). 
A function $\mu: \mathcal{Z}(H) \rightarrow [0,\infty]$ is called a \emph{cylindrical measure} on $H$, if for each finite subset $\Gamma \subseteq H$ the restriction of $\mu$ to $\mathcal{C}(H,\Gamma)$ is a measure. A cylindrical measure $\mu$ is said to be \emph{finite} if $\mu(H)< \infty$ and a \emph{cylindrical probability measure} if $\mu(H)=1$. 

A \emph{cylidrical random variable} over $H$ is a linear map $X: H \rightarrow L^{0} \ProbSpace$. We can associate to $X$ a cylindrical probability measure $\mu_{X}$ on $H$, called its \emph{cylindrical distribution}, and given by 
\begin{equation*} 
\mu_{X}(Z) \defeq \Prob \left( ( X(h_{1}), \dots, X(h_{n})) \in A  \right), 
\end{equation*}
% = \Prob \circ X^{-1} \circ \pi_{h_{1}, \dots, h_{n}}^{-1}(A)
for the cylindrical set $Z=\mathcal{Z}\left(h_{1}, \dots, h_{n}; A \right)$. 
Conversely, to every cylindrical probability measure $\mu$ on $H$ there is a canonical cylindrical random variable for which $\mu$ is its cylindrical distribution (see \cite{SchwartzRM}, p.256-8). 
%Analogous we say that a cylindrical random variable is $p$-\emph{integrable} if it is a linear map from $H$ to $L^p(\Omega,\mathbb{P})$. And we say that $X$ has \emph{zero-mean} if $ \Exp \left( X(h) \right)=0$,  $\forall \, h \in H$. 
Finally we call as the \emph{characteristic function} of a cylindrical random variable $X$ as the function $\phi:H\rightarrow\mathbb{C}$ defined as $\phi_X(h)=\mathbb{E}(\exp\{iX(h)\})$.

Let $J=[0,\infty[=\mathbb{R}_+$ or $J=[0,T]$ for some $T>0$. We say that $X=\{ X_{t}\}_ {t \in J}$ is a \emph{cylindrical process} if $X_{t}$ is a cylindrical random variable over $H$ for each $t \in J$ 
%(analogous a cylindrical process is $p$-\emph{integrable} if for every $t\in J$, $X_t$ is a $p$-integrable random variable)
. We say that $X$ has \emph{c\`adl\`ag paths} 
%or we say that is a  \emph{cylindrical square-integrable martingale} 
if for every $h\in H$ we have $\{X_t(h)\}_{t\in J}\in D_J(\mathbb{R})$.
%and $\{X_t(h)\}_{t\in J}\in\mathcal{M}_J^2(\mathbb{R})$ respectively.

To every $H$-valued stochastic process $X=\{X_{t}\}_{ t \in J}$ there corresponds a cylindrical process defined by $\{\inner{X_{t}}{h}\}_{ t \in J}$, for each $h \in H$. We will say that it is the \emph{cylindrical process induced} by $X$. Conversely, if a $H$-valued processes $Y=\{Y_{t}\}_{t \in J}$ is said to be a $H$-valued \emph{version} of the cylindrical process $X=\{X_{t}\}_{t \in J}$ if for each $t \in J$ and $h \in H$, $\inner{Y_{t}}{h}=X_{t}(h)$ $\Prob$-a.e.

%\begin{rema}
%The theorem \ref{LevyKhintCylTheo} in \cite{Applebaum-Riedle-2010} is proved for the context of cylindrical measures, but by the fact that a cylindrical process has a cylindrical law that it is itself a cylindrical measure also using the result mentioned in \cite{RiedleWeinerCyl-2018}, we can extended as we presented as before.
%\end{rema}

\section{Radonification of a cylindrical L\'evy process and applications}

\begin{defi}(see \cite{Applebaum-Riedle-2010}, p.705)\label{LevyCil}
Let $L=\{L_t\}_{t\in\mathbb{R}_+}$ a cylindrical process. We say that $L$ is a \emph{cylindrical L\'evy process} if for any $h_1,...,h_n\in H$ the process $\{(L_t(h_1),...,L_t(h_n))\}_{t\in\mathbb{R}_+}$ is a $\mathbb{R}^n$-valued L\'evy process.
\end{defi}

The following result is about the construction of a $H$-valued process that is a version of a cylindrical process. This construction has been used in more general context instead of a Hilbert space, so the same one is strongly inspired in \cite{Fonseca-Mora-RegularizationNuclear} with some differences in some arguments. However, because this result is vital in our followings results we present a short-proof, also we considered that because the context in \cite{Fonseca-Mora-RegularizationNuclear} is more general and so, more complicated it is worth it to present the proof in a Hilbert space context. Some of these tricks in the following calculations are well known as Gaussian bounds.

\begin{theo}[\textbf{Radonification}]\label{Radonification}
Let $X=\{X_t\}_{t\in\mathbb{R}_+}$ a cylindrical process with c\`adl\`ag paths (respectively continuous paths) such that for every $t\in\mathbb{R}_+$ the linear map $X_t:H\rightarrow L^0(\Omega,\mathcal{F},\mathbb{P})$ is continuous. If $S\in\mathcal{L}_2(H,G)$ then there exist a unique (up to indistinguishable versions) $G$-valued process $Y=\{Y_t\}_{t\in\mathbb{R}_+}\in D(G)$ (respectively $C(G)$) such that is a version of the cylindrical process $X\circ S^*=\{X_t\circ S^*\}_{t\in\mathbb{R}_+}$.
\end{theo}
\begin{prf}
We will prove it for \cadlag paths, the argument is the same for continuous paths. First we prove the statement for an interval $[0,T]$. We claim that $X:H\rightarrow D_T(\mathbb{R})$ defined by $h\mapsto\{X_t(h)\}_{t\in[0,T]}$ is linear continuous. In effect the map $X$ is well-defined and obviously linear. For the continuity we apply the closed graph theorem then we just need to prove that it is closed. Let $h_n\rightarrow h$ such that $X(h_n)\rightarrow Z$ in $D_T(\mathbb{R})$. Then, for each $t \in [0,T]$, $X_{t}(h_{n}) \rightarrow Z_{t}$ in probability. But since $h_{n} \rightarrow h$, we also have $X_{t}(h_{n}) \rightarrow X_{t}(h)$ in probability. By uniqueness of limits $Z_{t}= X_{t}(h)$ $\Prob$-a.e. Since both processes are c\`{a}dl\`{a}g, then $Z$ and $X(h)$ are indistinguishable and so then $X$ is closed and hence continuous.

With this continuity it is not difficult to prove that for every $\delta>0$ and every $h\in H$ there exist a $\rho>0$ such that:
\begin{equation}\label{boundRadon}
    \mathbb{E}\left(\displaystyle\sup_{t\in [0,T]}\left\vert 1-e^{iX_t(h)}\right\vert \right)\leq \delta +\frac{2 \Vert h \Vert ^2}{\rho^2}.
\end{equation}

Now we consider the decomposition of $S^*$ as in \eqref{spectra} and so we take the sequence in $G$ given by $\{Y_n\}_{n\in\mathbb{N}}:=\{\sum_{k=0}^n{\lambda_kX_t(h_k)}g_k\}_{n\in\mathbb{N}}$ for each $t\in[0,T]$. Clearly it is also a sequence in $D_T(G)$. We claim that is a Cauchy sequence in $D_T(G)$. Let $\epsilon>0$ and $n\geq m>0$, now first by Parseval's identity (using the complete orthonormal set $\{g_n\}_{n\in\mathbb{N}}$):
$$\mathbb{P}\left(\frac{1}{\epsilon}\sup_{t\in[0,T]}{\left\Vert Y_n-Y_m\right\Vert}>1\right)= \mathbb{P}\left(\frac{1}{\epsilon^2}\sup_{t\in[0,T]}{\sum_{k=m}^n{\vert \lambda_kX_t(h_k)g_k\vert^2}}>1\right).$$
 Now by Chebyshev's inequality and defining a $\mathbb{R}^{n-m}$-valued random variable with normal distribution  $Z\sim N(0,\epsilon^2I)$, also applying Tonelli's Theorem we obtain:
$$\mathbb{P}\left(\frac{1}{\epsilon^2}\sup_{t\in[0,T]}{\sum_{k=m}^n{\vert \lambda_kX_t(h_k)g_k\vert^2}}>1\right)\leq \frac{\sqrt{e}}{\sqrt{e}-1}\int_{\mathbb{R}^{n-m}}{\mathbb{E}\left(\sup_{t\in [0,T]}\vert 1-e^{iX_t\left(\sum_{k=n}^m{z_k\lambda_k h_k}\right)} \vert\right)}p_Z(d\Vec{z}).$$
Finally applying \ref{boundRadon} and Cauchy-Schwartz inequality we get:
$$\frac{\sqrt{e}}{\sqrt{e}-1}\int_{\mathbb{R}^{n-m}}{\mathbb{E}\left(\sup_{t\in [0,T]}\vert 1-e^{iX_t\left(\sum_{k=n}^m{z_k\lambda_k h_k}\right)} \vert\right)}p_Z(d\Vec{z})\leq \frac{\sqrt{e}}{\sqrt{e}-1}\left(\delta +\frac{2}{\rho^2 \epsilon^2}\sum_{k=m}^n{\vert \lambda_k \vert ^2}\right)$$
and by the square sumability of $\{\lambda_k\}_{k\in\mathbb{N}}$ and because $\delta$ is arbitrary we conclude by taking $n,m\rightarrow\infty$ that $\{Y_n\}_{n\in\mathbb{N}}$ is a Cauchy sequence in $D_T(G)$ and by completeness we obtain the existence and uniqueness of the $G$-valued process $Y=\{Y_t\}_{t\in[0,T]}$ defined for every $t\in[0,T]$ as $Y_t=\sum_{k=0}^\infty{\lambda_kX_t(h_k)g_k}$. To prove that $Y$ is a version of $X\circ S^*$ is a direct calculation using \eqref{spectra}. Finally to extend to $t\in\mathbb{R}_+$ just applying the result over $[0,T_n]$ for $T_n\rightarrow\infty$ and by the uniqueness of the $G$-valued process constructed it can be extended.
\end{prf}

\begin{lemm}\label{versionLevy}
Let $X=\{X_t\}_{t\in\mathbb{R}^+}$ a cylindrical L\'evy process over $H$. Then there exist a cylindrical L\'evy process $\tilde{X}=\{\tilde{X}_t\}_{t\in\mathbb{R}^+}$ with \cadlag paths over $H$ such that $X_t(h)=\tilde{X}_t(h)\mathrm{ \ \ }\mathbb{P}-a.s$ for all $h\in H$.
\end{lemm}
\begin{prf}
Let $h\in H$ so there exist a L\'evy process with \cadlag paths $\{\tilde{X}_t(h)\}_{t\in\mathbb{R}^+}$ such that is a version of $\{X_t(h)\}_{t\in\mathbb{R}^+}$. Then define $\tilde{X}$ as $\tilde{X}(h)=\{\tilde{X}_t(h)\}_{t\in\mathbb{R}^+}$. It is easy to verify that $\tilde{X}$ is a cylindrical L\'evy process.

Because $H$ is separable then for a orthonormal base $\{h_k\}_{k\in\mathbb{N}}$ we have:
\begin{equation*}
    X_t(h)=\tilde{X}_t(h) \Leftrightarrow \sum_{k=0}^\infty{\langle h_k,h\rangle X_t(h_k)}=\sum_{k=0}^\infty{\langle h_k,h\rangle \tilde{X}_t(h_k)}\Leftrightarrow X_t(h_k)=\tilde{X}_t(h_k)
\end{equation*}
so to prove that $X_t(h)=\tilde{X}_t(h)\mathrm{ \ \ }\mathbb{P}-a.s$ for all $h\in H$ it is enough to prove it over $\{h_k\}_{k\in\mathbb{N}}$.

Define $A_k:=\{\omega\in\Omega :X_t(h_k)=\tilde{X}_t(h_k)\}$ with  $\mathbb{P}(A_k)=1$. Finally, note that: 
$$\mathbb{P}\left(\bigcap_{k=0}^\infty{A_k}\right)=1-\sum_{k=0}^\infty{\mathbb{P}(A_k^c)}=1,$$
so $\mathbb{P}(X_t(h_k)=\tilde{X}_t(h_k),\forall k\in\mathbb{N})=1$.
\end{prf}

Now our main result:

\begin{theo}\label{radonLevy}
Let $X=\{X_t\}_{t\in\mathbb{R}_+}$ a cylindrical L\'evy process such that for every $t\in\mathbb{R}_+$ the linear map $X_t:H\rightarrow L^0(\Omega,\mathcal{F},\mathbb{P})$ is continuous. If $S\in\mathcal{L}_2(H,G)$ then there exist a unique (up to indistinguishable versions) $G$-valued L\'evy process $Y=\{Y_t\}_{t\in\mathbb{R}_+}\in D(G)$ (respectively $C(G)$ if the cylindrical L\'evy process has continuous paths) such that is a version of the cylindrical process $X\circ S^*=\{X_t\circ S^*\}_{t\in\mathbb{R}_+}$.
\end{theo}
\begin{prf}
We will prove the statement without the hypothesis of continuous paths because it is analogous. Let $X=\{X_t\}_{t\in\mathbb{R}_+}$ a cylindrical L\'evy process, then by lemma \ref{versionLevy} there exist a $\tilde{X}=\{\tilde{X}_t\}_{t\in\mathbb{R}_+}$ a cylindrical L\'evy process with \cadlag paths. Now by theorem \ref{Radonification} there exist a a unique $G$-valued process $Y=\{Y_t\}_{t\in\mathbb{R}_+}\in D(G)$ such that is a version of the cylindrical process $\tilde{X}\circ S^*=\{\tilde{X}_t\circ S^*\}_{t\in\mathbb{R}_+}$. So we just need to verify that $Y$ satisfy the definition of a L\'evy process.
\begin{itemize}
    \item Let $g\in G$, by linearity we have $0=\tilde{X}_0\circ S^*(g)=\langle Y_0,g\rangle$ then $Y_0=0$ almost surely.
    
    \item Let $0\leq t_1\leq \cdot\cdot\cdot \leq t_n$, as $\tilde{X}$ is a cylindrical L\'evy process then it has stationary increments for every $h\in H$. So for $g\in H$ we have:  
    $$\tilde{X}_{t_{j+1}}\circ S^*(g)-\tilde{X}_{t_j}\circ S^*(g)=(\tilde{X}_{t_{j+1}}-\tilde{X}_{t_j})\circ S^*(g)\stackrel{d}{=} \tilde{X}_{t_{j+1}-t_j}\circ S^*(g),$$
    for $j=1,...,n$. Then  $\phi_{\tilde{X}_{t_{j+1}}\circ S^*-\tilde{X}_{t_j}\circ S^*}(g)=\phi_{\tilde{X}_{t_{j+1}-t_j}\circ S^*}(g)$. Finally:
    \begin{eqnarray*}
    \phi_{Y_{t_{j+1}}-Y_{t_j}}(g)
    & = & \mathbb{E}\left(e^{i\langle Y_{t_{j+1}}-Y_{t_j},g \rangle}\right) = \mathbb{E}\left(e^{i(\tilde{X}_{t_{j+1}}\circ S^*-\tilde{X}_{t_j}\circ S^*)(g)}\right)=\phi_{\tilde{X}_{t_{j+1}}\circ S^*-\tilde{X}_{t_j}\circ S^*}(g)\\[1mm]
    & = & \phi_{\tilde{X}_{t_{j+1}-t_j}\circ S^*}(g)=\mathbb{E}\left(e^{i\tilde{X}_{t_{j+1}-t_j}\circ S^*(g)} \right)=\mathbb{E}\left(e^{i\langle Y_{t_{j+1}-t_j},g \rangle} \right)\\[1mm]
    & = & \phi_{Y_{t_{j+1}-t_j}}(g).
    \end{eqnarray*}
    Then we have $Y_{t_{j+1}}-Y_{t_j}\stackrel{d}{=}Y_{t_{j+1}-t_j}$.
    
    \item Let $0\leq t_1\leq \cdot\cdot\cdot \leq t_n$ y $j,m\in\{1,...,n\}$ with $j>m$, we must prove that  $Y_{t_{j+1}}-Y_{t_j}$ is independent of  $Y_{t_{m+1}}-Y_{t_m}$, that is  $\sigma((Y_{t_{j+1}}-Y_{t_j})^{-1}(A):A\in\mathcal{B}(G))$ is independent of  $\sigma((Y_{t_{m+1}}-Y_{t_m})^{-1}(A):A\in\mathcal{B}(G))$.
    
    Let $Z_g(B_1)$ and $Z_g(B_2)$ in $\mathcal{Z}(G)$, we will denote them by $Z_1$ y $Z_2$ to simplify. Then:
    \begin{align*}
    &\mathbb{P}((Y_{t_{j+1}}-Y_{t_j})^{-1}(Z_1)\cap (Y_{t_{m+1}}-Y_{t_m})^{-1}(Z_2))\\[1mm]
    &=\mathbb{P}(\{\langle Y_{t_{j+1}}-Y_{t_j},g \rangle \in B_1\} \cap \{\langle Y_{t_{m+1}}-Y_{t_m},g \rangle \in B_2\})\\[1mm]
    &=\mathbb{P}(\{ (\tilde{X}_{t_{j+1}}-\tilde{X}_{t_j})\circ S^*(g) \in B_1\} \cap \{(\tilde{X}_{t_{m+1}}-\tilde{X}_{t_m})\circ S^*(g) \in B_2\})\\[1mm]
    &=\mathbb{P}((\tilde{X}_{t_{j+1}}-\tilde{X}_{t_j})\circ S^*(g) \in B_1)\mathbb{P}((\tilde{X}_{t_{m+1}}-\tilde{X}_{t_m})\circ S^*(g) \in B_2)\\[1mm]
    &=\mathbb{P}((Y_{t_{j+1}}-Y_{t_j})^{-1}(Z_1))\mathbb{P}((Y_{t_{m+1}}-Y_{t_m})^{-1}(Z_2)).
    \end{align*}
    
    It is easy to verify the same calculation for an arbitrary choice of cylinder sets of any size. Then $Y_{t_{j+1}}-Y_{t_j}$ is independent of  $Y_{t_{m+1}}-Y_{t_m}$ over $\mathcal{Z}(G)$. Let  $\mathcal{A}_1:=\{(Y_{t_{j+1}}-Y_{t_j})^{-1}(Z):Z\in\mathcal{Z}(G)\}$ and  $\mathcal{A}_2:=\{(Y_{t_{m+1}}-Y_{t_m})^{-1}(Z):Z\in\mathcal{Z}(G)\}$. Here  $\mathcal{A}_1,\mathcal{A}_2\subset\mathcal{F}$, and how  $\mathcal{Z}(G)$ is an algebra, it follows that  $\mathcal{A}_1,\mathcal{A}_2$ is a $\pi$-system, then  $\sigma(\mathcal{A}_1)$ and $\sigma(\mathcal{A}_2)$ are  independient. We conclude that  $\sigma(\mathcal{Z}(G))=\mathcal{C}(G)=\mathcal{B}(G)$.
    
    \item Let $\epsilon>0$. For the stationary increments we note that for $t\geq s\geq 0$ we have:
    $$\mathbb{P}(\Vert Y_t-Y_s \Vert>\epsilon)=\mathbb{P}(\Vert Y_{t-s} \Vert>\epsilon),$$
    we called $l=t-s$.
    If $t\rightarrow s^+$ or $t\rightarrow s^-$ implies that $l\rightarrow 0^+$ and because $Y\in D(G)$ it follows that $\lim_{l\rightarrow 0^+}{Y_l}=Y_0=0$ and then:
    $$\lim_{t\rightarrow s^+}{\mathbb{P}(\Vert Y_t-Y_s \Vert>\epsilon)}=\lim_{t\rightarrow s^+}{\mathbb{P}(\Vert Y_{t-s} \Vert>\epsilon)}=0,$$
    and:
    $$\lim_{t\rightarrow s^-}{\mathbb{P}(\Vert Y_t-Y_s \Vert>\epsilon)}=\lim_{t\rightarrow s^-}{\mathbb{P}(\Vert Y_s-Y_t \Vert>\epsilon)}=\lim_{t\rightarrow s^-}{\mathbb{P}(\Vert Y_{s-t} \Vert>\epsilon)}=0.$$
\end{itemize}
\end{prf}

\begin{exam}(see \cite{Applebaum-Riedle-2010}, p.705)\label{WeinerCil}
Let $W=\{W_t\}_{t\in\mathbb{R}_+}$ a cylindrical process. We will say that $W$ is a \emph{cylindrical Weiner process} if for any $h_1,...,h_n\in H$ the process $\{(W_t(h_1),...,W_t(h_n))\}_{t\in\mathbb{R}_+}$ is a $\mathbb{R}^n$-valued Weiner process. We say that $W$ is a standard cylindrical Weiner process if every $\{(W_t(h_1),...,W_t(h_n))\}_{t\in\mathbb{R}_+}$ is a $\mathbb{R}^n$-valued standard Weiner process.
\end{exam}

\begin{exam}(see \cite{Applebaum-Riedle-2010}, p.706)\label{CompoundCil}
Let $Z=\{Z_n\}_{n\in\mathbb{N}}$ a cylindrical random variable sequence with same cylindrical distribution such that for each $h\in H$, $\{Z_n(h)\}_{n\in\mathbb{N}}$ are independent. Let $N=\{N_t\}_{t\in\mathbb{R}_+}$ a Poisson process with intensity $\lambda>0$, independent of $\{Z_n(h)\}_{n\in\mathbb{N}}$ for every $h\in H$.  We say that $X=\{X_t\}_{t\in\mathbb{R}_+}$ is a  \emph{cylindrical compound Poisson process} if for each $t\in\mathbb{R}_+$ it has the form of:
$$X_t= \left\{ \begin{array}{lcc}
             0 &   if  & N_t=0 \\
             \\ \displaystyle\sum_{k=0}^{N_t}Z_k &    & \mathrm{in \ other \ case} 
             \end{array}
   \right.$$
\end{exam}

\begin{coro}
Let  $W=\{W_t\}_{t\in\mathbb{R}^+}$ a cylindrical Wiener process such that for every $t\in\mathbb{R}_+$ the linear map $W_t:H\rightarrow L^0(\Omega,\mathcal{F},\mathbb{P})$ is continuous. If $S\in\mathcal{L}_2(H,G)$ then there exist a unique (up to indistinguishable versions) $G$-valued Wiener process $Y=\{Y_t\}_{t\in\mathbb{R}_+}\in C(G)$ such that is a version of the cylindrical process $W\circ S^*=\{W_t\circ S^*\}_{t\in\mathbb{R}_+}$.
\end{coro}
\begin{prf}
As we know $W$ is a L\'evy process with continuous paths so by theorem \ref{radonLevy} the statement is almost finished because $Y$ satisfies definition of a L\'evy process. We just need to prove that $Y$ have Gaussian increments but it follows immediately because for any $g\in G$ we have $\langle Y_t-Y_s,g\rangle=W_t(S^*(g))-W_s(S^*(g))$ and because $W\circ S^*(g)$ has Gaussian increments then $Y_t-Y_s$ has normal distribution.
\end{prf}

\begin{coro}
Let $X=\{X_t\}_{t\in\mathbb{R}^+}$ a cylindrical compound Poisson process such that for every $t\in\mathbb{R}_+$ the linear map $X_t:H\rightarrow L^0(\Omega,\mathcal{F},\mathbb{P})$ is continuous. If $S\in\mathcal{L}_2(H,G)$ then there exist a unique (up to indistinguishable versions) $G$-valued compound Poisson process $Y=\{Y_t\}_{t\in\mathbb{R}_+}\in D(G)$ such that is a version of the cylindrical process $X\circ S^*=\{X_t\circ S^*\}_{t\in\mathbb{R}_+}$.
\end{coro}
\begin{prf}
Because $X$ is a cylindrical compound Poisson process there exist a sequence $Z=\{Z_n\}_{n\in\mathbb{N}}$ of cylindrical random variable with same cylindrical distribution such that for each $h\in H$, $\{Z_n(h)\}_{n\in\mathbb{N}}$ are independent and $N=\{N_t\}_{t\in\mathbb{R}_+}$ a Poisson process with intensity $\lambda>0$, independent of $\{Z_n(h)\}_{n\in\mathbb{N}}$ for every $h\in H$ which satisfies the descomposition:
$$X_t= \left\{ \begin{array}{lcc}
             0 &   if  & N_t=0 \\
             \\ \displaystyle\sum_{k=0}^{N_t}Z_k &    & \mathrm{in \ other \ case} 
             \end{array}
   \right.$$
   Applying the operator $S^*$ on both sides for a $g\in G$:
   $$ X_t\circ S^*(g)= \left\{ \begin{array}{lcc}
             0 &   if  & N_t=0 \\
             \\ \displaystyle\sum_{k=0}^{N_t}Z_k\circ S^*(g) &    & \mathrm{in \ other \ case.} 
             \end{array}
   \right.$$
   Now applying theorem \ref{radonLevy} to $X$ and theorem $\ref{Radonification}$ to $Z$ it follows the existence of a $G$-valued L\'evy process $Y$ which is a version of $X\circ S^*$ and $G$-valued sequence of random variables $M=\{M_k\}_{k\in\mathbb{N}}$ which satisfies $Z_k\circ S^*(g)=\langle M_k,g\rangle$ such that:
   $$Y_t= \left\{ \begin{array}{lcc}
             0 &   if  & N_t=0 \\
             \\ \displaystyle\sum_{k=0}^{N_t}M_k &    & \mathrm{in \ other \ case} 
             \end{array}
   \right.$$
   To finished this result we must verify that $M$ satisfies the conditions of definition of a Compound Poisson process. To check that every random variable of $M$ has same distribution is immediately. The independence of every member of $M$ follows by a same argument of independence used in the proof of theorem \ref{radonLevy}. We just need to verify that $N$ and $M$ are independent processes.
   
   Let $0\leq t_1\leq\cdot\cdot\cdot\leq t_m$ and let  $\{A_n\}_{n\in\mathbb{N}}=\{A_g(B_n)\}_{n\in\mathbb{N}}$ cylinder sets. Let $j\in \{1,...,m\}$, by independence of  $N$ with $Z\circ S^*$ for every $h\in H$ we have:
   \begin{eqnarray*}
   \mathbb{P}\left(\{N_{t_j}=k\}\cap\left(\bigcap_{n\in\mathbb{N}}{M_n^{-1}(A_n)}\right)\right)
   & = & \mathbb{P}\left(\{N_{t_j}=k\}\cap\left(\bigcap_{n\in\mathbb{N}}\{\langle M_n,g\rangle\in B_n\}\right)\right)\\[1mm]
   & = & \mathbb{P}\left(\{N_{t_j}=k\}\cap\left(\bigcap_{n\in\mathbb{N}}\{Z_n\circ S^*(g)\in B_n\}\right)\right)\\[1mm]
   & = & \mathbb{P}\left(\{N_{t_j}=k\}\cap\left(\bigcap_{n\in\mathbb{N}}{(Z_n\circ S^*(g))^{-1}(B_n)}\right)\right)\\[1mm]
   & = & \mathbb{P}(N_{t_j}=k)\prod_{n\in\mathbb{N}}{\mathbb{P}{((Z_n\circ S^*(g))^{-1}(B_n))}}\\[1mm]
   & = & \mathbb{P}(N_{t_j}=k)\prod_{n\in\mathbb{N}}{\mathbb{P}{(M_n^{-1}(A_n))}}.
   \end{eqnarray*} 
   This calculation extends to an arbitrary choice of times for $N$ and any choice of cylinder sets of any size. The results is concluded thanks to the fact that $\mathcal{C}(G)=\mathcal{B}(G)$.
\end{prf}

\textbf{Acknowledgements} { The author acknowledge The University of Costa Rica for providing financial support through the grant of the thesis ``Procesos de L\'evy cil\'indrico: Regularizaci\'on e Integraci\'on''}.

\end{document}